\documentclass[11pt]{article}

\usepackage{a4wide}
\usepackage{amsmath}
\usepackage{amsfonts}
\usepackage{amssymb}
\usepackage{euscript}

\newtheorem{theorem}{Theorem}

\newtheorem{lemma}[theorem]{Lemma}
\newtheorem{corollary}[theorem]{Corollary}

\newenvironment{proof}{\noindent {\it Proof:~~}\ }{\  \rule{1mm}{2mm}\medskip}
\newenvironment{proofof}[2]{\noindent {\it Proof of #1}~#2: \ }{~\rule{1mm}{2mm}\medskip}

\def\transp#1{\,{^t\hspace{-1pt}#1}}
\def\esp#1{{\rm E}\left[#1\right]}
\def\prob#1{{\rm P}\!\left(#1\right)}

\def\vol#1{{\rm Vol}\left(#1\right)}

\newcommand{\Z}{\mathbb Z}
\newcommand{\F}{\mathbb F}
\newcommand{\R}{\mathbb R}
\renewcommand{\H}{\mathbf H}
\newcommand{\A}{\mathbf A}
\newcommand{\I}{\mathbf I}
\newcommand{\rand}{C_{\rm rand}}

\newcommand{\x}{\mathbf x}
\newcommand{\y}{\mathbf y}

\newcommand{\p}{\mathfrak p}
\renewcommand{\a}{\mathbf a}
\renewcommand{\u}{\mathbf u}

\begin{document}
\title{ On the construction of dense lattices with a given automorphism group}
\author{Philippe Gaborit and Gilles Z\'emor}
\date{May 23, 2006}

\maketitle
\begin{abstract}
We consider the problem of constructing dense lattices of $\R^n$ with
a given non trivial automorphism group.
We exhibit a family of such lattices
of density at least $cn2^{-n}$, which
matches, up to a multiplicative constant,
the best known density of a lattice packing.
For an infinite sequence of dimensions $n$, we exhibit a finite set
of lattices
that come with an automorphism group of size $n$, and
a constant proportion of which achieves the aforementioned
lower bound on the largest packing density.
The algorithmic complexity for exhibiting a basis of such a lattice
is of order $exp(n \log n)$, which improves upon previous theorems
that yield an equivalent lattice packing density.
The method developed here
involves applying Leech and Sloane's Construction~A
to a special class of codes with a given automorphism group, namely
the class of double circulant codes.
\end{abstract}


\section{Introduction}

A lattice packing of Euclidean balls in $\R^n$
is a family of disjoint Euclidean balls of equal radius 
centered on the points of some non degenerate lattice. 
The proportion of the space covered by these Euclidean balls
is called the density of the packing. 
When balls of volume $V$ are packed by a lattice
$\Lambda$, the corresponding density is $V.\det(\Lambda)^{-1}$, where 
$\det(\Lambda)$ denotes the determinant of the lattice, i.e. the volume
of a fundamental region.

The classical Minkowski-Hlawka Theorem  states that for 
$n$ greater than $1$ there exist lattice packings with density
at least $2^{1-n}\zeta(n)$. 
This lower bound on the lattice packing density was later improved
by a linear
factor to a quantity of the form $cn2^{-n}$ for constant~$c$.
This improvement is originally due to Rogers \cite{rog47} with $c=2e^{-1}$.
The constant~$c$ was successively improved by
Davenport and Rogers \cite{dr47} to $c=1.68$ and eventually by Ball
\cite{bal92} to $c=2$.

In the meantime, Rush \cite{rus89}, building upon a technique of Rush and Sloane
\cite{rs87}, essentially recovered the original Minkowski-Hlawka lower bound
on the largest density of a sphere packing
using coding theory arguments together with the Leech-Sloane
Construction A for lattices. While this did not achieve the
improved density of the form $cn2^{-n}$, it had the alternative advantage
of being more effective than the proofs of the above results.
Rush's construction
exhibits in a natural way
a finite number of lattices among which dense ones exist. 
This number, though still too large to be in any way practical, is
much smaller than what can be derived by applying the original proofs of the
results highlighted above: consequently, the algorithmic complexity of
Rush's construction is
of the form $exp(n\log n)$ which is a substantial improvement over
the preceding ones (see~\cite{cs88}, p.18).

Recently, the $cn2^n$ improved lower bound on the minimum density was
made as effective as Rush's lattice construction, with $c=0.01$, 
for (non-lattice) sphere packings 
by Krivelevich, Litsyn and Vardy in \cite{klv04}.
They use an elegant graph theory method that enables them to find dense 
packings with a time (and space) complexity $exp(n \log n)$.

In this paper, we again make the $cn2^n$ lower bound as effective, with
$c \approx 0.06$, without
paying the price of losing lattice structure. In fact, the dense
lattice packings that we exhibit have {\em additional} algebraic structure,
namely they come together with an automorphism group of size $n$.
This additional structure is not a by-product of our method but is
an essential reason for the improved density. This is a small step
towards showing that, in the asymptotic setting,
 algebraic constructions can compete with
unstructuredness, and maybe even stand out.

The starting point of our approach is similar to that of \cite{rs87} 
and \cite{rus89}, namely relies upon construction A to transform
codes in $\F_p^n$ into lattices of~$\R^n$. The specificity of the
Rush-Sloane method is to consider codes designed for a metric
which is unconventional in $\F_p^n$ but
specially adapted to the Euclidean metric in~$\R^n$.
However, instead of indiscriminately looking for the best codes for this metric
in the whole space~$\F_p^n$, we depart from \cite{rs87,rus89}
by restricting our attention to an exponentially smaller 
set of codes, namely a class that has a 
given automorphism group (double circulant codes),
and prove that a constant fraction of them yield lattices
with improved density. Similar codes were also used in a coding theory
context to improve the classical
Gilbert-Varshamov bound for linear codes by a linear factor~\cite{gz06}.
Exhibiting a lattice basis has algorithmic (time) complexity $exp(n \log n)$.

The paper is organized as follows: in Section \ref{sec:dense}, we show how
dense lattices are constructed from ``dense'' codes and we formulate
our main results, Theorem~\ref{th:VGplus} and Corollary \ref{cor:lattice}.
In Section~\ref{sec:double}
we show how to obtain good double circulant codes.

\section{From dense codes to dense lattices}\label{sec:dense}
Let $S_n$ denote the Euclidean ball of radius $1$ in $\R^n$, we have:
\begin{equation}
  \label{eq:volume}
  \vol{S_n} = \frac{\pi^{(n/2)}}{(n/2)!}.
\end{equation}
Let $S_n(d)$ denote the Euclidean ball of radius $d$ in $\R^n$, so that we
have $$\vol{S_n(d)}=d^n\vol{S_n}.$$

Let $\rho\in\R$ be the radius of a
Euclidean ball of volume $p^{n/2}$ for $p$ any positive number,
i.e. $\vol{S_n(\rho)}=p^{n/2}$.
By \eqref{eq:volume} and Stirling's formula we have:
\begin{equation}
  \label{eq:radius}
  \rho = \sqrt{\frac{pn}{2e\pi}}(1+o(1))
\end{equation}
where $o(1)$ will always be understood to mean a quantity that vanishes
as $n$ goes to infinity.

For $\Lambda$ a lattice of dimension $n$ it is customary to define its minimum
norm by 
$$\mu(\Lambda)=\min\left\{\sum_{i=1}^nx_i^2 \: ,\; (x_1,\cdots,x_n) 
  \in \Lambda\setminus\{0\}\right\}.$$
The lattice $\Lambda$ defines a packing of $\R^n$ by spheres of
Euclidean radius $\sqrt{\mu}/2$ and the 
density of this packing is given by:
\begin{equation}
  \label{eq:Delta}
  \Delta=\frac{\vol{S_n(\sqrt{\mu}/2)}}{\det(\Lambda)}
         =\frac{\vol{S_n}\mu^{n/2}}{2^n\det(\Lambda)},
\end{equation}
where $\det(\Lambda)$ stands for the determinant of $\Lambda$.

From now on let $p$ be a prime. 
We identify elements $z$ of $\F_p$ with elements $z$ of $\Z$ such
that $$-\frac{p-1}{2}\leq z \leq \frac{p-1}{2}.$$
With this convention, following \cite{rs87,rus89}, 
we introduce the {\em norm} of
a vector $\x=(x_1,\cdots,x_n)$ in $\F_p^n$ as the non-negative real number
$$\|\x\|_2=\sqrt{x_1^2+\cdots +x_n^2}.$$

Let $B_{n,p}(d)$ denote the set of vectors $\x\in\F_p^n$ such that
$\|\x\|_2\leq d.$
We shall only be dealing with values of $d$ such that
$d<p/2$ so that we shall always have:
\begin{equation}
  \label{eq:ball}
  |B_{n,p}(d)| = |\Z^n\cap S_n(d)|
\end{equation}
hence, by fitting the sphere $S_n(d)$ inside a union of $n$-cubes of volume $1$,
\begin{equation}
  \label{eq:size}
  \vol{S_n\left(d-\frac{\sqrt{n}}{2}\right)} \leq 
 |B_{n,p}(d)| \leq \vol{S_n\left(d+\frac{\sqrt{n}}{2}\right)}.
\end{equation}
Let us call a $[n,k,d,p]$ code a $k$-dimensional subspace $C$
of $\F_p^n$ such that $d$ equals the minimum of the norm $\|\x\|_2$ 
of a nonzero codevector $\x\in C$. We will refer to $d$ as the
minimum norm of the code $C$.

Recall that Construction A associates to a code $C$ the lattice:
$$ A_p(C)=\{(x_1, ...,x_n) \in \Z^n \: | \: 
            (x_1 \: \bmod \: p, ...,x_n \: \bmod \: p) 
\in C\}.$$ 
It is readily seen that this lattice has
minimum norm $\mu=min(d^2,p^2)$ and determinant $p^{n-k}$.
In the following, we will always ensure that
$d\leq p$ so that the $[n,k,d,p]$ code $C$ yields by construction A a lattice
of $\R^n$ of norm $d^2$ with density~\eqref{eq:Delta}:
\begin{equation}
  \label{eq:density}
  \Delta = \frac{1}{2^n}\frac{\vol{S_n(d)}}{p^{n-k}}.
\end{equation}
By \eqref{eq:size} this gives a density
\begin{equation}
  \label{eq:density2}
  \Delta \geq \frac{1}{2^n}\frac{|B_{n,p}(d)|}{p^{n-k}}
                           \left(1+\frac{\sqrt{n}}{2d}\right)^{-n}.
\end{equation}
We shall prove
\begin{theorem}\label{th:VGplus}
  There exists a constant $c$, such that for any $n=2q$, $q$ a large
enough prime, there exists  a prime $p$, $n^2\log n<p\leq (n^2\log^2 n)^{5.5}$,
and an $[n,n/2,d,p]$ code $C$ such that
  $$|B_{n,p}(d)| \geq cnp^{n/2}.$$
Furthermore, the automorphism group of $C$ contains a subgroup isomorphic
to $\Z/2\Z\times \Z/q\Z$.
\end{theorem}

The condition $n^2\log n<p$ in Theorem \ref{th:VGplus} will ensure that
the term $(1+\sqrt{n}/2d)^{-n}$ in \eqref{eq:density2}
tends to $1$ when $n$ tends to infinity.
This will enable us to obtain:

\begin{corollary}\label{cor:lattice}
There exists a constant $c$, such that for any $n=2q$, $q$ a large
enough prime, there exists  a
lattice of $\R^n$ with density at least $cn/2^n$
and whose automorphism group contains a subgroup isomorphic
to $\Z/2\Z\times \Z/q\Z$. Such a lattice can be constructed 
with time complexity $exp(n \log n))$.
\end{corollary}

The numerical value of the constant $c$ in Theorem \ref{th:VGplus}
and Corollary \ref{cor:lattice} can be estimated to be at least
$(2-1/e)(2+e^2\pi)^{-1}\approx 0.064.$

\section{Double circulant codes and random choice}\label{sec:double}
A $p$-ary {\em double circulant code} is a $[2q,q,d,p]$ linear code $C$ with a 
parity-check matrix of the form $\H = [\I_q\, |\,\A]$ where $\I_q$ is the
$q\times q$ identity matrix and 
$$\A = 
\begin{bmatrix}
       a_1 & a_2 & \dots & a_q\\
       a_q & a_1 & \dots & a_{q-1}\\
       a_{q-1} & a_q & \dots & a_{q-2}\\
       \hdotsfor{4} \\
       a_2 & a_3 & \dots & a_1
\end{bmatrix}.$$
This simply means that $C$ is the kernel of the mapping 
$\x\mapsto \x\transp{\H}$ from $\F_p^{2q}$ to $\F_p^q$.
  
We will only consider the case when $q$ is a prime. Let $n=2q$.
There is a natural action of the group $G=\Z/2\Z\times\Z/q\Z$ on the space
$\F_p^{n}$ of vectors $\x=(x_1\ldots x_q,x_{q+1}\ldots x_{2q})$
\begin{eqnarray*}
    G\times \F_p^n & \rightarrow & \F_p^n\\
           (g,\x)     & \mapsto     & g\cdot\x
\end{eqnarray*}
where $(0,1)\cdot\x =(x_q,x_1\ldots x_{q-1},x_{2q},x_{q+1},\ldots x_{2q-1})$
and $(1,0)\cdot\x = -\x$.
The double circulant code $C$ is invariant under this group action and
so is the norm of any vector $\x$.
Note that construction A applied to the code $C$ will clearly yield a lattice
whose automorphism group contains $G$.

To show that double circulant codes with a large minimum norm $d$ exist,
we shall study the typical behaviour of $d$ when a double circulant
code is chosen at random. We now formalize this:

Consider the random double circulant code $\rand$ obtained by choosing the
the first row of $\A$, the vector $(a_1\ldots a_q)$, 
with a uniform distribution in $\F_p^q$.
We are interested in the random variable $X(w)$ equal to the number
of nonzero codevectors of $\rand$ of norm not more than $w$. In other
words we define
  $$X(w) = \sum_{\x\in B_{n,p}(w)\setminus\{0\}}X_\x$$
where $X_\x$ is the Bernoulli random variable
equal to $1$ if $\x\in \rand$ and equal to zero otherwise.
Our strategy is to study
the maximum value of $w$ for which we can claim
$\prob{X(w)>0}<1$, this will prove the existence of codes of parameters
$[n,n/2,d> w,p]$. 

The core remark is now that, if $\y=g\cdot\x$, then
  $$X_\y = X_\x.$$
Let now $B_{n,p}'(w)$ be a set of representatives of the orbits of the
elements of $B_{n,p}(w)$, i.e. for any $\x\in B_{n,p}(w)$, 
$|\{g\cdot \x ,g\in G\}\cap B_{n,p}'(w)|=1$. We clearly have $X(w)>0$
if and only if $X'(w)>0$ where
  $$X'(w) = \sum_{\x\in B_{n,p}'(w)\setminus\{0\}}X_\x.$$
Denote by $\ell(\x)$ the length (size) of the orbit of $\x$,
i.e. $\ell(\x)=\#\{g\cdot \x ,g\in G\}$. We have
\begin{equation}
  \label{eq:orbits}
  X'(w) = \sum_{\x\in B_{n,p}(w)\setminus\{0\}}\frac{X_\x}{\ell(\x)}
\end{equation}
By writing $\prob{X(w)>0}=\prob{X'(w)>0}\leq \esp{X'(w)}$, together
with \eqref{eq:orbits} we obtain
\begin{equation}
  \label{eq:orbits2}
  \prob{X(w)>0}\leq\sum_{\lambda|n}
  \sum_{\substack{\|\x\|_2\leq w\\ \ell(\x)=\lambda}}
  \frac{\esp{X_\x}}{\lambda}.
\end{equation}
Since $n=2q=|G|$ and $q$ is a prime, possible values of $\lambda$ in
\eqref{eq:orbits2} are $1,2,q,n$. Note that $\ell(0) = 1$,
$\ell(\x)=2$ for $\x$ of the form 
$\x =(\alpha(1,1,\ldots 1),\beta(1,1,\ldots 1))$ and $\ell(\x)\geq q$ for
all other vectors. In fact a closer look shows that
$\ell(\x) = q$ is not possible. For this to happen, one of the two
halves of $\x$, call it $\y$,
would have all its $q$ cyclic shifts distinct, and the property that
$-\y$ equals some cyclic shift of $\y$. But then it would be possible to
partition the set of cyclic shifts of $\y$ into pairs of opposite vectors,
but $q$ is odd, a contradiction.
Therefore Inequality \eqref{eq:orbits2} gets rewritten as
\begin{equation}
  \label{eq:orbits3}
  \prob{X(w)>0}\leq
  \sum_{\substack{\x =(\alpha(1,1,\ldots 1),\beta(1,1,\ldots 1))\\
  0<\|\x\|_2\leq w}}
  \frac{\esp{X_\x}}{2} +
  \sum_{\substack{\ell(\x)=n\\ 0<\|\x\|_2\leq w}}
  \frac{\esp{X_\x}}{n}.
\end{equation}
We now switch to evaluating the right hand side of \eqref{eq:orbits3}.

\subsection{Syndrome distribution}
We need to study carefully the quantities
$\esp{X_\x}=\prob{\x\in\rand}$, for $\x\in B_{n,p}(w)$.
For $\x\in \F_p^{n}$, let us write $\x = (\x_L, \x_R)$ with
$\x_L,\x_R\in \F_p^q$. Consider the syndrome function $\sigma$
\begin{eqnarray*}
 \sigma : \F_p^{n} & \rightarrow & \F_p^q\\
          \x     & \mapsto     & \sigma(\x) =\x\transp{\H}
                                 =\sigma_L(\x)+\sigma_R(\x)
\end{eqnarray*}
where $\sigma_L(\x)=\x_L$ and $\sigma_R(\x)=\x_R\transp{\A}$.

For any vector $\u=(u_0,\ldots ,u_{q-1})$ of $\F_p^q$, 
denote by $\u(Z)=u_0 + u_1Z+\cdots +u_{q-1}Z^{q-1}$ its polynomial 
representation in the ring ${\mathbf R}=\F_p[Z]/(Z^q-1)$. 
For any $\u\in \F_p^q$,
let $C(\u)$ denote the cyclic code of length $q$ generated by the
polynomial representation of $\u$ (i.e. $C(\u)$ is the ideal generated
by $\u(Z)$ in the ring ${\mathbf R}$). We have:
\begin{lemma}\label{lem:Cx}
The right syndrome  $\sigma_R(\x)$ of any given $\x\in \F_p^{n}$
is uniformly distributed in the cyclic code $C(\x_R)$. Therefore, the
probability $\prob{\x\in\rand}$ that $\x$ is a codevector of the random
code $\rand$ is
$$\begin{array}{ll}
\bullet\hspace{2mm} \prob{\x\in\rand} = 1/|C(\x_R)| & 
                    \text{if}\hspace{2mm} \x_L\in C(\x_R),\\
\bullet\hspace{2mm}  \prob{\x\in\rand} = 0       & 
                    \text{if}\hspace{2mm} \x_L\not\in C(\x_R).
\end{array}$$
\end{lemma}

\begin{proof}
A little thought shows that $\sigma_R(\x)$ has polynomial representation
equal to $\x_R(Z)\a(Z)$, where $\a=(a_1,a_q,a_{q-1},\ldots ,a_2)$ is the 
transpose of the first column of $\A$. Therefore,  the image of the mapping
\begin{eqnarray*}
\psi~:\F_p^q & \rightarrow & \F_p^q\\
  \a     & \mapsto     & \sigma_R(\x)
\end{eqnarray*}
for fixed $\x$, is the cyclic code $C(\x_R)$. Since this mapping is linear,
every element of $C(\x_R)$ has the same number of preimages (namely
$Ker\psi$), therefore when the distribution of $\a$ is uniform in
$\F_p^q$, the distribution of $\sigma_R(\x)$ is uniform in the code $C(\x_R)$.
\end{proof}

\subsection{The choice of $p$ and the cyclic codes $C(\x_R)$}\label{sec:p}
The right hand side of \eqref{eq:orbits3} will be easiest to study if there
are as few as possible cyclic codes in $\F_p^q$, i.e. if the ring
${\mathbf R}$ has as few as possible invertible elements, equivalently
if $Z^q-1$ has as few as possible divisors in $\F_p[Z]$.
The next lemma tells us how to ensure this, while simultaneously bounding
from above the size of~$p$, so as to retain some control over the overall
contruction complexity.
\begin{lemma}\label{lem:p}
  For any $n=2q$ large enough, there exists a prime $p$ in the range 
$n^2\log n\leq p \leq (n^2\log^2 n)^{5.5}$ for which the
the factorization of $Z^q-1$ into irreducible polynomials of $\F_p[Z]$ is
  $$Z^q-1 = (Z-1)(1+Z+Z^2+\cdots +Z^{q-1}).$$
\end{lemma}
\begin{proof}
  We just need to find $p$ in the required range such that
  $(p \bmod q)$ is a primitive element in $\Z/q\Z$.

Let $Q=q^2\p$ where 
$\p$ is a prime such that $4\log n\leq \p\leq 4\log^2 n$:
$\p$ exists for $q$ large enough, and we have
$n^2\log n\leq Q\leq n^2\log^2n$.

Let $\alpha<q$ be a positive integer that is a primitive element in $\Z/q\Z$.
Since $q$ is prime we have $q\neq 0\bmod \p$ so that we may choose
$\varepsilon_1\in\{1,2\}$ and $\varepsilon_2\in\{0,1\}$ such that
$r = (1+\varepsilon_1 q)(\alpha +\varepsilon_2q)$ is
coprime to $\p$ and therefore to $Q$. Note also that $r$
is smaller than $Q$ for $q$
large enough, not prime, and equal to $\alpha \bmod q$.
By Linnik's
Theorem on least primes in arithmetic progressions, there
exists a prime $p$ such that $p=r \;\bmod Q$ and $p\leq Q^L$
for a constant~$L$.
We have $p=r=\alpha \;\bmod q$. Note that since $r$ is not
prime we have $Q<p$ in addition to $p\leq Q^L$.
By a result of Heath-Brown \cite{hb92} we have $L\leq 5.5$.
\end{proof}

For $p$ as in Lemma \ref{lem:p}
we therefore have exactly two non-trivial cyclic codes over $\F_p$
of length $q$,
namely $C_1$, the subspace generated by the all-one vector (or
the generator polynomial $1+Z+\cdots +Z^{q-1}$) and its dual, $C_1^\perp$,
with generator polynomial $Z-1$.

Now Lemma \ref{lem:Cx} implies that there are exactly two types of
non-zero vectors of $\F_p^n$ such that $\prob{\x\in\rand}$ is different
from zero and from $1/p^q$, namely:
\begin{itemize}
\item vectors $\x$ such that $\x_L\in C_1$ and $\x_R\in C_1$, we call them
      vectors of type $1$. For these vectors we have $\prob{\x\in\rand}=1/p$.
\item vectors $\x$ such that $\x_L\in C_1^\perp$ and $\x_R\in C_1^\perp$, 
      we call them vectors of type $2$.
      For these vectors we have $\prob{\x\in\rand}=1/p^{q-1}$.
\end{itemize}
Next, we study the number of these exceptional vectors to evaluate
their contribution to the upper bound \eqref{eq:orbits3}.
\subsection{Number of vectors of type 1 and type 2 in $B_{n,p}(\rho)$}
Suppose $w=\rho(1+o(1))$ where $\rho$ is defined by $S_n(\rho)=p^q$
(see \eqref{eq:radius}). Note that in Lemma \ref{lem:p} we have
chosen $p$ such that \eqref{eq:radius} implies
$\sqrt{n}/\rho = o(1/n)$. Therefore,~\eqref{eq:size} implies in turn that
\begin{equation}
  \label{eq:equivalent}
  |B_{n,p}(w)|=\vol{S_n(w)}(1+o(1)).
\end{equation}
A vector of type~$1$ in $B_{n,p}(w)$ is a vector $\x$ such that
  $$\x_L = \alpha(1,1,\ldots ,1)\hspace{5mm}\text{and}\hspace{5mm}
    \x_R = \beta(1,1,\ldots ,1).$$
The number $N_1(w)$
of possible values of $(\alpha,\beta)$ such that $\|\x\|_2\leq w$
is, 
  $$N_1(w)=
        \#\{(\alpha,\beta)\in\F_p^2\;|\;\alpha^2\frac n2 +\beta^2\frac n2\leq 
        w^2\}.$$
Therefore, for $w<(p-1)/2$ (which is always going to be satisfied for
$n$ large enough and $p$ chosen as in Lemma \ref{lem:p}),
  $$N_1(w)=\#\{(\alpha,\beta)\in\Z^2\;|\;\alpha^2+\beta^2
               \leq\frac{2w^2}{n}\}$$
and, bounding from above by the area of a $2$-dimensional disc,
  $$N_1(w)\leq \pi\left(w\sqrt{\frac 2n}+\sqrt{2}\right)^2.$$
Therefore \eqref{eq:radius} gives
\begin{equation}
  \label{eq:N1}
  N_1(w) \leq \frac pe (1+o(1)).
\end{equation}

We now switch to evaluating the cardinality $N_2(w)$ of the set
$A$ of vectors of type $2$ in $B_{n,p}(w)$.
Now let $B$ be the set of vectors $\y$ of $\F_p^n$ obtained by the following
procedure:
\begin{enumerate}
\item choose $\x=(x_1\ldots x_n)\in A$
\item choose $i,j$ with $1\leq i\leq q$, $q+1\leq j\leq 2q$
\item choose two integers $l,r$ such that $|l|\leq \lceil\sqrt{tp}\rceil$ and
      $|r|\leq \lceil\sqrt{tp}\rceil$, where $t$ is a constant to be determined
      later
\item define $\y=(y_1\ldots y_n)$ by $y_i=l$, $y_j=r$ and $y_h=x_h$ for 
      $h\neq i,j$.
\end{enumerate}
We now define the bipartite graph with vertex set $A\cup B$ by putting
an edge between $\x\in A$ and $\y\in B$ if $\y$ is obtained from
$B$ by the above procedure. Let $E$ be the set of edges of this graph.
The degree of a vertex $\x\in A$ is clearly 
$q^2(2\lceil\sqrt{tp}\rceil+1)^2\geq 4tpq^2$ so that we have
$|E|\geq |A|4tpq^2$. Recall that $\x$ is of type $2$ means that
$x_1+\cdots +x_q=0$ and $x_{q+1}+\cdots +x_{2q}=0$.
Now let $\y\in B$. There is at most one way of modifying
two given coordinates $i,j$, $1\leq i\leq q$, $q+1\leq j\leq 2q$,
so as to obtain a vector $\x\in A$. In other words the degree of a
vertex $\y\in B$ is at most $q^2$ and $|E|\leq |B|q^2$.
We have therefore
\begin{equation}
  \label{eq:4p|B|}
  |A|\leq \frac{1}{4tp} |B|.
\end{equation}
Now notice that if $\x\in A$ and $\y\in B$ are adjacent in the bipartite
graph we have
  $$\|\y\|_2^2 \leq \|\x\|_2^2+2\lceil\sqrt{tp}\rceil^2$$
so that $B\subset B_{n,p}(w')$ with 
$w'=\sqrt{w^2+2\lceil\sqrt{tp}\rceil^2}$.
Since $w=\rho(1+o(1))$, this gives
\begin{equation}
  \label{eq:w'}
  w' = w\sqrt{1+2tp\rho^{-2}(1+o(1))}.
\end{equation}
In particular we have $w'=\rho(1+o(1))$ so that,
applying \eqref{eq:equivalent}, we get
\begin{align*}
  |B|\leq |B_{n,p}(w')| & = \vol{S_n(w')}(1+o(1)) = \vol{S_n}w'^n(1+o(1))\\
     & = \vol{S_n(w)}\frac{w'^n}{w^n}(1+o(1))
       =|B_{n,p}(w)|\frac{w'^n}{w^n}(1+o(1)).
\end{align*}
Now \eqref{eq:w'} and \eqref{eq:radius} give:
$$|B|\leq |B_{n,p}(w)|\left(1+\frac{4te\pi}{n}\right)^{n/2}(1+o(1)).$$
Together with \eqref{eq:4p|B|} we obtain the following bound on
$N_2(w)=|A|$:
$$  N_2(w)\leq \frac{e^{2te\pi}}{4tp}|B_{n,p}(w)|(1+o(1)).$$
Now choose $t=(2e\pi)^{-1}$ so as to minimize $e^{2te\pi}/4tp$ and we get:
 \begin{equation}
  \label{eq:N2}
  N_2(w)\leq \frac{e^{2}\pi}{2p}|B_{n,p}(w)|(1+o(1)).
\end{equation}

\subsection{Proof of Theorem \ref{th:VGplus} and Corollary \ref{cor:lattice}}
We are now ready to prove the main result.

\medskip

\begin{proofof}{Theorem}{\ref{th:VGplus}}
Choose $p$ as in Lemma \ref{lem:p} and
choose $w$ such that $\vol{S_n(w)} = cnp^{n/2}$, $c$ a constant to be
determined later. This clearly implies $w=\rho(1+o(1))$ so that,
by \eqref{eq:equivalent}, we have $|B_{n,p}(w)|=cnp^{n/2}(1+o(1))$.
The upper bounds \eqref{eq:N1} and \eqref{eq:N2} apply
and~\eqref{eq:orbits3} yields:
\begin{align*}
  \prob{X(w)>0}&\leq N_1(w)\frac{p^{-1}}{2}+N_2(w)\frac{p^{1-n/2}}{n}
      +|B_{n,p}(w)|\frac{p^{-n/2}}{n}\\
  &\leq \frac{1}{2e} + \frac{e^{2}\pi}{2}c + c +o(1).
\end{align*}
We obtain therefore $\prob{X(w)>0}<1$ for $n$ large enough and
  $$c < \frac{2-\frac{1}{e}}{2+e^2\pi}\approx 0.064.$$
We have proved that for such a value of $c$, some
double circulant codes with minimum norm $d\geq w$ must exist.
\end{proofof}

\begin{proofof}{Corollary}{\ref{cor:lattice}}
Let $C$ be the code in Theorem \ref{th:VGplus}.
By inequality \eqref{eq:size}, since $cnp^{n/2}\leq |B_{n,p}(d)|$,
the quantity $d+\sqrt{n}/2$ must be greater than the radius of a Euclidean
ball of volume $cnp^{n/2}$. As before,
by equality \eqref{eq:radius}, the code's
minimum norm $d$ must be greater than $\sqrt{pn}$ 
multiplied by a constant, so that the term
  $\left(1+\sqrt{n}/2d\right)^{-n}$
in \eqref{eq:density2} converges to $1$ when $n\rightarrow\infty$,
since $\sqrt{p}/n\rightarrow \infty$. Therefore~\eqref{eq:density2}
yields the announced density for the lattice deduced from the code~$C$
by construction~A. 

Construction~A preserves the automorphism group
of the code in the lattice. 
The construction complexity is simply that
of going over all double circulant codes
of length $n$ over $\F_p$ (there are $p^{\frac{n}{2}}$ of them), 
and checking, by exhaustive search over
the $p^{\frac{n}{2}}$ codevectors, 
whether they contain a vector of norm less than the
required bound. The resulting complexity equals therefore
$p^n$ times quantities of a lesser order of magnitude, i.e.
$p^{n(1+o(1))}$ which is not more, by Lemma \ref{lem:p},
than $2^{2L(1+o(1))n\log_2(n)}$.
\end{proofof}

\section{Concluding comments}

\begin{itemize}
\item 
The proof of Theorem \ref{th:VGplus} shows that, by
lowering the value of $c$, we can make all the contributions to
the probability of the existence of a codevector of weight $\leq w$
vanish, except for the codevectors of type $1$.
In other words, for small values of the constant~$c$, the asymptotic
probability that the double circulant code-random lattice yields a packing of density less 
than $cn2^{-n}$ equals the non-vanishing
probability (not more than $1-\frac{1}{2e}$) that codevectors of 
type~$1$ exist.
When this happens, not only does the packing density drop below
$cn2^{-n}$, but it drops below the Minkowski density altogether.
In contrast, typical random lattice packings have a density of order
$1/2^n$ \cite{st01}.
\item
The action of the automorphism group of the lattices presented here
is not transitive on the set of coordinates, it has two orbits.
Can one construct dense lattices with a transitive automorphism group~?
\item The automorphism group here has size (at least) $n$. Could alternative
constructions yield an automorphism group of guaranteed larger size (potentially
resulting in increased packing densities)~?
\end{itemize}

% \bigskip
% \noindent {\bf Philippe Gaborit}\\
% Universit\'e de Limoges\\
% XLIM-DMI\\
% 123, av. Albert Thomas,\\
% 87000 Limoges, France.\\
% mail:gaborit@unilim.fr\\

% \bigskip

% \noindent {\bf Gilles Z\'emor}\\
% \'Ecole Nationale Sup\'erieure des T\'el\'ecommunications\\
% CNRS UMR 5141,\\
% D\'ept Informatique et R\'eseaux,\\
% 46 rue Barrault,\\
% 75634 Paris Cedex 13, France.\\
% mail:zemor@enst.fr

\end{document}